\documentclass[12pt]{amsart}
\usepackage{amsmath}
\usepackage{amssymb,amsthm,mathrsfs}
\usepackage{commath} 
\usepackage{tikz-cd}
\usepackage{hyperref}
\usepackage{cleveref}

\newtheorem{thm}{Theorem}[section]
\newtheorem{cor}[thm]{Corollary}
\newtheorem{lem}[thm]{Lemma}
\newtheorem{prop}[thm]{Proposition}

\theoremstyle{definition}

\newtheorem{remark}[thm]{Remark}

\numberwithin{equation}{section}

\let\P\relax
\DeclareMathOperator{\P}{\mathbb{P}}

\DeclareMathOperator{\Z}{\mathbb{Z}}

\DeclareMathOperator{\cM}{\mathcal{M}}
\DeclareMathOperator{\cC}{\mathcal{C}}

\DeclareMathOperator{\Jac}{Jac}

\DeclareMathOperator{\coker}{coker}

\newcommand{\cls}[1]{\overline{ #1 }}
\newcommand{\iso}{\cong}
\DeclareMathOperator{\Hom}{Hom}

\DeclareMathOperator{\Hilb}{Hilb}

\DeclareMathOperator{\tensor}{\otimes}

\DeclareMathOperator{\End}{End}
\DeclareMathOperator{\tor}{tor}

\newcommand{\wt}[1]{\widetilde{#1}}

\DeclareMathOperator{\into}{\hookrightarrow}
\DeclareMathOperator{\dashto}{\dashrightarrow}

\title{Threefolds containing all curves are rationally connected}
\author{Sixuan Lou}

\address{Department of Mathematics, Statistics, and CS \\
University of Illinois at Chicago, Chicago IL 60607}
\email{slou2@uic.edu}

\begin{document}

\begin{abstract}
    Any smooth projective curve embeds into $\P^3$. More generally, any curve
    embeds into a rationally connected variety of dimension at least three. We
    prove conversely that if every curve embeds in a variety $X$, then $X$
    contains a rationally connected subvariety of dimension at least three. 
    In particular, ``all curves embed'' is a birational property for threefolds.
\end{abstract}

\maketitle

\section{Introduction}
Let $k$ be an algebraically closed field of characteristic zero.

A common theme in the study of algebraic varieties is to determine their
properties by studying the subvarieties contained in them. In particular, rational
curves play a central role in the birational geometry of projective varieties.
A variety $X$ is \emph{rationally connected} if for a general pair
of points of $X$ there is a rational curve passing through them.
The abundance of rational curves in these varieties allows curves in them to
deform in large families. Using this observation, Sankaran has shown that
rationally connected varieties contain all curves \cite{S11}.

In this paper, we study the converse question: what can we say about a variety
that admits a family of morphisms from curves with high variation?
Recall that the \emph{variation} of a family of smooth curves of genus $g$ is
the dimension of the image of the classifying map to the moduli stack
$\mathcal{M}_g$.
We will show that if the variation of curves satisfies a specific lower
bound, the variety must be uniruled.

\begin{thm}[Theorem \ref{main-thm-ext}]
    Let $X$ be a smooth projective variety of dimension $n > 1$, and let $g \geq 2$
    be an integer. Assume $X$ admits a covering family of morphisms from genus $g$
    curves with variation $\nu$, whose general member is birational onto its
    image. If $\nu > n + g - 2$, then $X$ is uniruled.
\end{thm}

Using this result as an essential ingredient, we obtain our main theorem, which
addresses the converse question for embeddings of general curves. We establish
that if a variety is swept out by embeddings of curves of sufficiently high
genus, it must contain a rationally connected subvariety of dimension at least
three.

\begin{thm}[Theorem \ref{general-RC}]
    Let $X$ be a smooth projective variety of dimension $n \geq 3$, and let $g > g(n)
    := \max(21, n - 2, (n+1)/2)$ be an integer. If $X$ admits a family of
    embeddings from genus $g$ curves that vary maximally in moduli, then $X$
    contains a rationally connected subvariety of dimension at least $3$.
\end{thm}

The proof of Theorem \ref{general-RC} proceeds by induction on the dimension of
the variety $X$. The crucial base case relies on establishing structural results
for surfaces that admit such families of morphisms.

\begin{thm}[Theorem \ref{surface}]
    Let $X$ be a smooth projective surface, and let $g \geq 2$ be an integer. If
    there is a covering family of morphisms from genus $g$ curves with maximal
    variation, whose general member is birational onto its image, then $X$ is
    rational.  Furthermore, if $g > 21$, then the general member of this family
    of morphisms cannot be an embedding.
\end{thm}

An immediate consequence of Theorem \ref{general-RC} relates to the
birational geometry of threefolds. By Theorem \ref{general-RC}, a smooth
projective threefold is rationally connected if and only if it admits an
embedding of every smooth projective curve. Since rational connectedness is a
birational invariant, it follows that the property of ``containing all curves''
must also be a birational property for smooth projective threefolds.

In general, if $\pi : \wt{X} \to X$ is a birational morphism between smooth
projective varieties, an embedding of a curve $C \into \wt{X}$ intersecting the
exceptional locus of $\pi$ does not necessarily descend to an embedding in $X$.
However, the global equivalence with rational connectedness allows us to bypass
this local obstruction entirely. If we can embed a general curve of sufficiently
high genus into $\wt{X}$, we can conclude that $\wt{X}$ is rationally connected.
This implies $X$ is rationally connected, and therefore, by Sankaran's result, $X$
must also contain all curves.

\begin{cor}
    Let $\pi : \wt{X} \to X$ be a birational morphism between smooth projective
    threefolds. If, for some $g > 21$, a general curve of genus $g$ embeds
    into $\wt{X}$, then both $\wt{X}$ and $X$ are rationally connected.
    Therefore, every curve has an embedding into $X$.
\end{cor}

We will first prove Theorem \ref{main-thm-ext} by examining the dimension of the
space of global sections of the normal sheaf $N_f$ of a general morphism. The
result will be established using the generalization of Clifford's theorem
bound to semistable bundles due to Brambila-Paz, Grzegorczyk and Newstead
\cite{BGN97}.  We then apply Theorem \ref{main-thm-ext} to the case of families of
surfaces and use it to study higher-dimensional cases.

\section{Proofs}

\begin{lem}
    \label{curve}
    Let $g \geq 2$ and let $M$ be a closed irreducible subvariety of the coarse
    moduli space $M_g$ of genus $g$ curves. If $\dim M \geq 2g - 1$, then a
    general member of $M$ does not cover any smooth projective curves other than
    $\P^1$ and itself.
\end{lem}

\begin{proof}
    By \cite{CGT92}, if $\dim M \geq 2g - 1$, then a general member $C$ of $M$
    has $\End \Jac(C) = \Z$; in particular, its Jacobian is simple.  Suppose $f
    : C \to D$ is a non-constant morphism to another smooth projective curve
    $D$. The morphism induces the pullback and pushforward morphisms
    \[
      f^* : \Jac(D) \to \Jac(C)  , \quad
      f_* : \Jac(C) \to \Jac(D)
    \]
    with $f_* \circ f^* = [\deg f]$. Thus, $\ker f^* \subseteq \Jac(D)[\deg f]$
    is finite.

    Let $A$ denote the image of $f^*$ in $\Jac(C)$. It is an abelian subvariety
    of $\Jac(C)$, isogenous to $\Jac(D)$.  Since $\Jac(C)$ is simple, either $A
    = 0$ or $A = \Jac(C)$. If $A = 0$, then $\Jac(D) = 0$ as well. Thus, $D \iso
    \P^1$. If $A = \Jac(C)$, then $C$ and $D$ have the same genus $g$.

    By the Riemann-Hurwitz formula,
    \[
        2g - 2 \geq (\deg f) (2g - g)
    \]
    Since $g \geq 2$, one concludes that $\deg f = 1$ and $f$ must be an
    isomorphism.
\end{proof}

\begin{thm}
    \label{main-thm-ext}
    Let $X$ be a smooth projective variety of dimension $n > 1$, and let $g \geq 2$
    be an integer. Assume $X$ admits a covering family of morphisms from genus $g$
    curves with variation $\nu$, whose general member is birational onto its
    image. If $\nu > n + g - 2$, then $X$ is uniruled.
\end{thm}

\begin{proof}
    Since $H_2(X, \Z)$ is countable, we may assume the family of morphisms maps
    to curves of a fixed class $\beta \in H_2(X, \Z)$.

    Let $\cC \to Z$ be a flat family of connected smooth projective curves of
    genus $g$ over a normal irreducible scheme $Z$. Let $F : \cC \to X \times Z$
    be the family of morphisms over $Z$. Let $F_X : \cC \to X$ be the projection
    to $X$, which is the universal evaluation map of the family of morphisms.

    Let $p : Z \to \cM_g(X, \beta)$ be the characteristic map to the moduli
    stack of stable maps. By \cite[Lemma 2.5]{DS17}, the scheme theoretic image
    $F(\cC)$ is flat over $Z$. Hence, we have a characteristic map $q : Z \to
    \Hilb$ to the Hilbert scheme of curves with an appropriate Hilbert
    polynomial in $X$.  Let $\pi : \cM_g(X, \beta) \to \cM_g$ be the projection.
    Since the family $\cC \to Z$ is equigeneric, and each morphism $F_t : C_t
    \to X$ is birational onto its image, the composition $\pi \circ p$ factors
    through a rational map $\varphi : q(Z) \dashto \cM_g$.
    \[
        \begin{tikzcd}
                      & Z \rar{p} \ar[ld,swap, "q"] \dar{q} & \mathcal{M}_g(X, \beta) \rar{\pi} & \mathcal{M}_g\\
            \Hilb  & q(Z) \ar[rru, dashed, swap, "\varphi"] \ar[l, swap, hook'] & &
        \end{tikzcd}
    \]

    By assumption, we have:
    \begin{enumerate}
        \item the image of $(\pi \circ p)(Z)$ has dimension $\nu$,
        \item the universal evaluation map $F_X : \cC \to X$ is dominant.
    \end{enumerate}

    Let $\eta \in Z$ be a general point.
    Let $N_{\eta}$ denote the normal sheaf of $F_{\eta}$. Let $N_{\eta, \tor}$
    denote the torsion subsheaf of $N_{\eta}$. Then we have the exact sequence
    \[
        0 \to N_{\eta, \tor} \to N_{\eta} \to \cls{N}_{\eta} \to 0.
    \]
    In particular, $\deg \cls{N}_{\eta} \leq \deg N_{\eta}$.

    Since $\dim (\pi \circ p)(Z) = \nu$, the differential
    \[
        d\pi_{p(\eta)} \circ dp_{\eta} : T_{\eta} Z \to H^0(C_{\eta}, N_{\eta})
        \to H^1(T_{C_{\eta}})
    \]
    has rank $\nu$.
    By \cite[Lemma 2.5]{DS17}, the preimage $dp_{\eta}^{-1}(H^0(N_{\eta, \tor}))
    $ lies in the kernel of $dq_{\eta}$. Hence,
    \[
        \nu \leq h^0(\cls{N}_{\eta}).
    \]

    We show $K_X \cdot \beta < 0$ if $\nu > n - 2 + g$. Assume otherwise
    that $K_X \cdot \beta \geq 0$. 

    By deformation theory, given an element $[f : C \to X, p]$ of $\cls{\cM}_{g,
    1}(X, \beta)$, the tangent space of $\cls{\cM}_{g,1}(X, \beta)$ at $[f,p]$
    is naturally identified with $H^0(\coker(df : T_C(-p) \to f^* T_X))$.
    Furthermore, the differential of the universal evaluation map $dF_X$ at $[f,
    p]$ is identified with the restriction map
    \[
        dF_X |_{(f,p)} : H^0(\coker(T_C(-p) \to
        f^*T_X)) \to \del{\coker(T_C(-p) \to f^* T_X)} \tensor k(p)
    \]

    By condition (2), the differential $dF_X |_{(\eta,p)}$ is surjective at a
    general point $\eta \in Z$ and a general point $p \in C_{\eta}$.
    This shows that the cokernel $Q := \coker(T_C(-p) \to F_{\eta}^* T_X) $ is
    generically globally generated. Since $\cls{N}_{\eta}$ is a quotient of $Q$,
    it is also generically globally generated.

    Consider the Harder-Narasimhan filtration of $\cls{N}_f$:
    \[
        0 = HN_0 \subsetneq HN_1 \subsetneq \cdots \subsetneq HN_{k-1}
        \subsetneq HN_{k} = \cls{N}_f
    \]
    Let $A_i := HN_i / HN_{i-1}$ for $i = 1, \ldots, k$. Let $r_i$ and $d_i$ denote
    the rank and the degree of each $A_i$, respectively, and let $\mu_i := d_i / r_i$
    denote its slope. We know each $A_i$ is semistable with
    \[
        \mu(A_1) > \mu(A_2) > \cdots > \mu(A_k).
    \]

    Since $A_k$ is a quotient of $\cls{N}_f$, it is generically globally
    generated. In particular, $A_k$ must admit a global section. Since $A_k$ is
    semistable, we conclude $\mu_k \geq 0$. Therefore, $\mu_i \geq 0$ for all $i
    = 1, \ldots, k$.

    Let $k_0$ be the integer between $0$ and $k$ such that
    \begin{align*}
        & \mu_i > 2g - 2, \quad \text{ for } 1 \leq i \leq k_0\\
        & \mu_i \leq 2g - 2, \quad \text{ for } k_0 + 1 \leq i \leq k.
    \end{align*}

    Since each $A_i$ is semistable, if its slope $\mu_i$ is greater than $2g - 2$,
    then by Serre duality we have $h^1(A_i) = \dim \Hom(A_i, \omega_C) = 0$.
    Moreover, by the extension of Clifford's theorem to semistable bundles on
    curves \cite[Theorem 2.1]{BGN97}, we obtain an upper bound for $h^0(A_i)$
    for $k_0 + 1 \leq i \leq k$, and hence an upper bound for $h^0(\cls{N}_f)$:
    \begin{align*}
        \nu & \leq h^0(\cls{N}_f) \leq \sum_{i = 1}^{k} h^0(A_i)
        = \sum_{i = 1}^{k_0} h^0(A_i) + \sum_{i = k_0 + 1}^{k} h^0(A_i)\\
        & \leq \sum_{i = 1}^{k_0} r_i (\mu_i + 1 - g) + \sum_{i = k_0 + 1}^{k}
        r_i + \frac{d_i}{2}
    \end{align*}

    If $k_0 = 0$, then the inequality simplifies to
    \[
        \nu \leq (n-1) + \frac{1}{2} \del{-K_X \cdot C + 2g - 2}
        \leq (n-1) + g - 1
    \]
    Hence, $\nu \leq n - 2 + g$.

    If $k_0 > 0$, we obtain the following:
    \[
        \nu
         \leq (n - 1) + \sum_{i = 1}^{k} d_i - g \sum_{i = 1}^{k_0} r_i
         \leq (n - 1) + 2g - 2 - g \sum_{i = 1}^{k_0} r_i
         \leq n - 3 + g
    \]

    Therefore, we must have $K_X \cdot \beta < 0$ if $\nu > n - 2 + g$.

    Assume $\nu > n - 2 + g$. Then $Z$ gives a covering family of irreducible
    curves of $X$ with negative intersection with $K_X$. We conclude that $K_X$
    is not pseudo-effective by \cite[Theorem 0.2]{BDPP13}. Therefore, $X$ is
    uniruled by \cite[Corollary 0.3]{BDPP13}.
\end{proof}

\begin{thm}
    \label{main-thm}
    Let $X$ be a smooth projective variety of dimension $n > 1$, and let $g$ be an
    integer with $g > (n+1)/2$. If $X$ admits a family of morphisms from genus
    $g$ curves with maximal variation, whose general member
    is birational onto its image, then $X$ is uniruled.
\end{thm}

\begin{proof}
    Apply \Cref{main-thm-ext} with $\nu = 3g - 3$. Therefore, $X$ must be
    uniruled if $g > (n+1)/2$.
\end{proof}

\begin{remark}
    \label{moduli-uniruled}
    Recall that the coarse moduli space of curves $\cls{M}_g$ is of general type
    for $g \geq 24$, and has Kodaira dimension $\geq 0$ for $g = 23$, by the
    foundational work of Harris, Mumford, and Eisenbud \cite{HM82, EH87}. More
    recently, Farkas, Jensen, and Payne improved this bound to show that $\cls{M}
    _g$ is of general type for all $g > 21$ \cite{FJP25}. Since a variety with
    non-negative Kodaira dimension cannot be uniruled, we conclude that the
    moduli space $\cls{M}_g$ is not uniruled for any $g > 21$.
\end{remark}

We now apply this general uniruledness criterion to establish the structural
result for surfaces.

\begin{cor}
    \label{surface}
    Let $X$ be a smooth projective surface, and let $g \geq 2$ be an integer. If
    there is a covering family of morphisms from genus $g$ curves with maximal
    variation, whose general member is birational onto its image, then $X$ is
    rational.  Furthermore, if $g > 21$, then the general member of this family
    of morphisms cannot be an embedding.
\end{cor}

\begin{proof}
    By Theorem \ref{main-thm}, one concludes $X$ is uniruled. Hence, $X$ is
    birationally ruled by some dominant rational map $p : X \dashto B$ to some
    smooth projective curve $B$.
    Since the curves in this family vary maximally in moduli, they must dominate $B$.
    Thus, from Lemma \ref{curve}, it follows $B$ must be rational, hence $X$
    is a rational surface. In this case, the family of morphisms must vary in a
    fixed linear system $V$. If the general member of this family of morphisms is an
    embedding, we obtain a dominant rational map $V \dashto M_g$.
    This
    contradicts the fact that the moduli space $\cls{M}_g$ is not uniruled for $g > 21$
    (see Remark \ref{moduli-uniruled}).
\end{proof}

This structural result of curves in surfaces extends naturally to families of
surfaces.

\begin{prop}
    \label{family-surface}
    Let $\pi : X \to B$ be a smooth morphism of relative dimension $2$ to a
    quasi-projective variety $B$ of dimension $b$. Let $g$ be an integer with $g
    > \max(21, b + 1)$.  Then there is no covering family of embeddings from
    genus $g$ curves that vary maximally in moduli to fibers of $\pi$.
\end{prop}

\begin{proof}
    Let $\cC \to Z$ be a family of embeddings from genus $g$ curves that cover
    $X$ and vary maximally in moduli. From the assumption, there is a morphism
    $\phi : Z \to B$ such that for each closed point $b \in B$, the restricted
    family over $Z_b = \phi^{-1}(b)$ is a covering family of embeddings from
    genus $g$ curves to the fiber $X_b = \pi^{-1}(b)$.

    Since $Z$ has maximal variation $3g - 3$, the fiber $Z_b$ has variation at
    least $3g - 3 - b$.  Since $g > b + 1$, we have $3g - 3 - b > g$. Therefore,
    by Theorem \ref{main-thm-ext}, each fiber $X_b$ is uniruled, hence must be
    ruled over some smooth curve $D_b$.  Since members of the family $Z_b$ are
    all embeddings, the family of morphisms $\cC_b \to Z_b$ must also dominate
    $D_b$. Moreover, since $Z_b$ has variation at least $3g - 3 - b \geq 2g - 1$,
    we conclude that $D_b$ is rational by Lemma \ref{curve}. Hence, $X_b$ is a
    rational surface for every $b \in B$. Therefore, the family of morphisms in
    $Z_b$ must vary in a fixed linear system $V_b$ on $X_b$. This gives a
    dominant rational map $V_b \dashto Z_b$. Since each fiber $Z_b$ is uniruled,
    the entire family $Z$ is also uniruled. Since $Z$ dominates $M_g$, this
    implies that the moduli space $\cls{M}_g$ is uniruled.
    We derive a
    contradiction when $g > 21$ (see Remark \ref{moduli-uniruled}).
\end{proof}

With these structural results established, we are now ready to
proceed with the induction to prove the main theorem for varieties of dimension
three and higher.

\begin{thm}
    \label{general-RC}
    Let $X$ be a smooth projective variety of dimension $n \geq 3$, and let $g > g(n)
    := \max(21, n - 2, (n+1)/2)$ be an integer. If $X$ admits a family of
    embeddings from genus $g$ curves that vary maximally in moduli, then $X$
    contains a rationally connected subvariety of dimension at least $3$.
\end{thm}

\begin{proof}
    We induct on the dimension of $X$. Let $\cC \to Z$ denote the family of
    embeddings. Let $F : \cC \to X$ denote the evaluation map. First observe
    that we may assume the family of embeddings sweeps out $X$. Otherwise, it must
    sweep out some subvariety $X' \subsetneq X$ of lower dimension. After
    resolving singularities of $X'$, we have a covering family of embeddings of
    genus $g$ curves to $\wt{X'}$ that vary maximally in moduli. Since $g(n)
    \geq g(n-1)$, by the inductive hypothesis we may reduce to the base case
    where $\wt{X'}$ is a surface. However, when $g > 21$, this is impossible by
    Proposition \ref{surface}.

    Let $\rho : X \dashto R(X)$ be the MRC fibration. The MRC quotient $R(X)$
    is not uniruled by \cite[Corollary 1.4]{GHS03}. Observe that if $\dim R(X) \leq n
    - 3$, then a very general fiber of $\rho$ will be an irreducible rationally
    connected subvariety of $X$, which establishes the conclusion. Therefore, we
    only need to consider the cases $\dim R(X) \in \{n, n - 1, n-2\}$.

    \begin{enumerate}
        \item If $\dim R(X) = n$, then $\rho$ is birational; hence, $X$ is not
            uniruled. We get a contradiction using Theorem \ref{main-thm} given
            $g > (n + 1) / 2$.

        \item Assume $\dim R(X) = n - 1$. Since the fibers of $\rho$ form a
            family of curves with bounded dimension $\leq n - 1 < 3g - 3$,
            the general member of $Z$ cannot be contracted by $\rho$.  Let $\eta \in
            Z$ be a general member, and let $F_{\eta} : C_{\eta} \to X$ be the
            corresponding embedding. Since $R(X)$ is projective, the rational
            map $\rho \circ F_{\eta} : C_{\eta} \to X \dashto R(X)$ can be
            completed to a morphism. Let $D_{\eta} \subseteq R(X)$ denote the
            image of $\rho \circ F_{\eta}$ and let $\wt{D}_{\eta}$ denote its
            normalization.  Then we get a covering of curves $C_{\eta} \to \wt{D}
            _{\eta}$.

            Since $Z$ dominates $\cM_g$, the curve $C_{\eta}$ is a general curve
            of $\cM_g$. Therefore, the curve $\wt{D}_{\eta}$ is either $\P^1$ or
            isomorphic to $C_{\eta}$ by Lemma \ref{curve}. Assume $D_{\eta} = \P^1$
            for a general member $\eta \in Z$. Since the composition $\rho \circ
            F : \cC \dashto R(X)$ is dominant, this would imply $R(X)$ is
            uniruled, which is a contradiction. Therefore, we conclude $C_{\eta}
            \to R(X)$ is birational onto its image for a general member $\eta \in Z$.

            After possibly resolving singularities of $R(X)$, we get a covering
            family of morphisms from genus $g$ curves that vary maximally in
            moduli whose general member is birational to its image in $R(X)$.
            Hence, from Theorem \ref{main-thm}, we derive the contradiction that
            $R(X)$ is uniruled once $g \geq n / 2$.

        \item Assume $\dim R(X) = n - 2$. Then $\cC \to Z$ induces a covering family
            of morphisms from genus $g$ curves varying maximally in moduli to
            the resolution $\wt{R(X)}$. By the same argument as the previous
            case, a general member $C_{\eta} \subseteq \cC$ must map
            birationally into $\wt{R(X)}$. Thus, we get a covering family of
            morphisms from genus $g$ curves with maximal variation, whose
            general member is birational to its image in $\wt{R(X)}$.
            Hence, from Theorem \ref{main-thm}, we derive the
            contradiction that $R(X)$ is uniruled once $g > (n - 1) / 2$.

            Thus, a general member of $Z$ is contracted by $\rho$.
            Let $R(X)^{\circ} \subseteq R(X)$ be the open dense subvariety such
            that the restriction $\rho^{\circ} : X^{\circ} = \rho^{-1}(R(X)
            ^{\circ}) \to R(X)^{\circ}$ is a smooth morphism. After restricting
            $Z$ to an open set, we obtain a covering family of embeddings of
            genus $g$ curves that vary maximally in moduli to fibers of
            $\rho^{\circ}$. We get a contradiction from Proposition \ref{family-surface}
            once $g > \max(21, n - 2)$.
    \end{enumerate}
\end{proof}

\begin{cor}
    Let $X$ be a smooth projective threefold. Then $X$ is rationally connected
    if and only if every smooth projective curve admits an embedding into $X$.
\end{cor}

\section{Competing interests}
No competing interest is declared.

\section{Acknowledgments}

The author would like to thank his advisor, \.{I}zzet Co\c{s}kun, for suggesting
this problem, for his constant patience, and for many helpful conversations
throughout the work. The author is also grateful to the referees for their
careful reading of the paper, numerous corrections that helped improve the
exposition, and specifically for recommending the extension of the technique
to prove the general result in Theorem \ref{main-thm-ext}.

\end{document}